\documentclass[12pt]{amsart}

\usepackage{amsmath}
\usepackage{amssymb}
\usepackage{amsthm}
\usepackage{amscd}
\usepackage{amsfonts, euscript,graphics}

\def \RM{\mathbb{R}}
\def \NM{\mathbb{N}}
\def \KM{\mathbb{K}}
\def\longermapsto{{\mapstochar\relbar\mskip-6mu\rightarrow}}
\def \.{\mskip 1mu}
\def \?{\mskip -1mu}

\newcommand{\ind}{\mathop{\rm ind}}

\newtheorem{thm}{Theorem}
\newtheorem{lem}{Lemma}
\newtheorem{prop}{Proposition}
\newtheorem{defi}{Definition}
\newtheorem{cor}{Corollary}

\theoremstyle{remark}
\newtheorem{rem}{Remark}
\newtheorem{example}{Example}

\title{Positive isotopies of Legendrian submanifolds and applications}

\author{Vincent COLIN, Emmanuel FERRAND, Petya PUSHKAR}

\date{}

\begin{document}

\begin{abstract}
We show that there is no positive loop inside the component of a fiber in
the space of Legendrian embeddings in the contact manifold $ST^*M$,
provided that the universal cover of $M$ is $\RM^n$. 
We consider some related results in the space of one-jets of functions on a compact manifold. 
We give an application to the positive isotopies in homogeneous neighborhoods of surfaces in a tight contact 3-manifold.
\end{abstract}

\maketitle

\section{Introduction and formulation of the results}

\subsection{} \label{intro}On the Euclidean unit $2$-sphere, the set of points which are at a
given distance of the north pole is in general a circle. When the
distance is $\pi$, this circle becomes trivial: it is reduced to the
south pole. Such a focusing phenomenon
cannot appear on a surface of constant, non-positive curvature.
In this case the image by the exponential map of a unit circle of vectors tangent to the
surface at a given point is never reduced to one point.

\medskip

In this paper, we generalize this remark in the context of contact
topology\footnote{in particular, no Riemannian structure is involved}.
Our motivation comes from the theory of the orderability of the group of
contactomorphisms of Eliashberg, Kim and Polterovitch \cite{EKP}.

\medskip

\subsection{Positive isotopies} Consider a $(2n+1)$-dimensional
manifold $V$ endowed with a {\em cooriented} contact structure $\xi$.
At each point of $V$, the contact hyperplane then separates the tangent
space in {\em a positive and a negative side}.

\begin{defi} \label{isoPos}
A smooth path $L_t=\varphi_t(L) , t\in [0,1]$ in the space of Legendrian embeddings
(resp. immersions) of a $n$-dimensional compact manifold $L$ in $(V,\xi)$ is called {\em a
  Legendrian isotopy (resp. homotopy)}.
If, in addition, for every $x \in L$ and every $t \in [0,1]$, the
velocity vector $\dot{\varphi_t}(x)$ lies in the {\em positive} side
of $\xi$ at $\varphi_t(x)$, then this Legendrian isotopy (resp. homotopy)
will be called {\em positive}.
\end{defi}

\medskip

\begin{rem}
This notion of positivity does not depend on the parametrization of
the $L_t$'s.
\end{rem}

If the cooriented contact structure $\xi$ is induced by a globaly defined contact form $\alpha$, the
above condition can be rephrased as $\alpha(\dot{\varphi_t}(x))>0$.
In particular, a positive contact hamiltonian induces positive isotopies.
 
\medskip
A positive isotopy (resp. homotopy) will be also called a {\em
positive path} in the space of Legendrian embeddings (resp. immersions).

\begin{example} \label{example1}

The space $J^1(N)=T^*N \times \RM $ of one-jets of functions on a
$n$-dimensional manifold $N$ has a natural contact one form
$\alpha=du-\lambda$, where $\lambda$ is the Liouville one-form of
$T^*N$ and $u$ is the $\RM$-coordinate. The corresponding contact
structure will be denoted by $\zeta$. Given a smooth function 
$f \colon N \to \RM$, its one-jet extension $j^1f$ is a Legendrian
submanifold. A path between two functions gives rise to an isotopy of
Legendrian embeddings between their one-jets extensions.

\medskip

A path $f_{t, t \in [0,1]}$ of functions on $N$ such that,
for any fixed $q \in N$, $f_t(q)$ is an increasing function of $t$,
gives rise to a positive Legendrian isotopy $j^1f_{t, t \in [0,1]}$ in $J^1(N)$.

\medskip

Conversely, one can check that a positive isotopy consisting only of
one-jets extensions of functions is always of the above type.
{\em In particular there are no positive loops consisting only of one-jet
extensions of functions}.

\end{example}

\begin{example} \label{ex1}
Consider a Riemannian manifold $(N,g)$.
Its unit tangent bundle $\pi\colon S_1N \to N$
has a natural contact one-form: If $u$ is a unit tangent vector to $N$,
and $v$ a vector tangent to $S_1N$ at $u$, then
$$\alpha (u) \cdot v=g(u, D\pi (u) \cdot v).$$
The corresponding contact structure will be denoted by $\zeta_1$.
The constant contact Hamiltonian $h=1$ induces the geodesic flow.

\medskip

Any fiber of $\pi\colon S_1N \to N$ is Legendrian. 
Moving a fiber by the geodesic flow is a typical example of a positive path.

\end{example}

\subsection{Formulation of results}
\label{results}
Let $N$ be a closed manifold.

\begin{thm} \label{thm1}
There is no closed positive path in the component of the space of
Legendrian embeddings in $(J^1(N), \zeta)$ containing the one-jet extensions of functions.
\end{thm}

The Liouville one-form of $T^*N$ induces a contact distribution of the
fiber-wise spherization $ST^*N$. This contact structure is
contactomorphic to the $\zeta_1$ of Example \ref{ex1}.  Our generalization
of the introductory Remark \ref{intro} is as follows.

\begin{thm} \label{thm2}
There is no positive path of Legendrian embeddings between two
distinct fibers of 
$\pi \colon ST^*N \to N$, provided that the universal cover of $N$ is $\RM^n$.
\end{thm}

\begin{thm} \label{thm3}
{\bf 0).} Any compact Legendrian submanifold of $J^1(\RM^n)$ belongs to a closed path of Legendrian embeddings.

{\bf i).} There exists a component of the space of Legendrian embeddings
in $(J^1(S^1), \zeta)$ whose elements are homotopic to $j^10$ and which contains a closed positive path.

{\bf ii).} There exists a closed positive path in the component of the space of
Legendrian immersions in $(J^1(S^1), \zeta)$ which contains the one-jet
extensions of functions.

{\bf iii).} Given any connected surface $N$, there exists a positive path of
  Legendrian immersions between any two fibers of $\pi\colon ST^*N \to N$.
\end{thm}

\medskip

F. Laudenbach \cite{La} proved recently the following generalization
of Theorem \ref{thm3} {\bf ii)}: for any closed $N$, there exists a
closed positive path in the component of the space of
Legendrian immersions in $(J^1(N), \zeta)$ which contains the one-jet
extensions of functions.

\medskip

Theorem \ref{thm3} {\bf 0)} implies that for any contact manifold $(V,\xi)$, there exists a closed positive path
of Legendrian embeddings (just consider a Darboux ball and embed the example of Theorem \ref{thm3} {\bf 0}).


\medskip

A Legendrian manifold $L \subset (J^1(N),\zeta)$ will be called {\em positive} 
if it is connected by a positive path to the one-jet extension of the zero function. 
The one-jet extension of a positive function is a positive Legendrian manifold. 
But, in general, the value of the $u$ coordinate can be negative at some points of a positive Legendrian manifold.

\medskip

Consider a closed manifold $N$ and fix a function $f\colon N \to {\RM}$. 
Assume that $0$ is a regular value of $f$.
Denote by $\Lambda$ the union for $\lambda \in {\RM}$ of the $j^1(\lambda f)$. 
It is a {\em smooth embedding} of ${\RM} \times N$ in  $J^1(N)$, foliated by the $j^1(\lambda f)$.
We denote by $\Lambda_+$ the subset $\bigcup_{\lambda>0}(j^1(\lambda f)) \subset \Lambda$.

\medskip

Consider the manifold $M=f^{-1}([0,+\infty[) \subset N$. Its boundary $\partial M$ is the set $f^{-1}(0)$. 
Fix some field $\KM$ and denote by $b(f)$ the total dimension of the homology of $M$ with
coefficients in that field ($b(f)=\dim_{\KM} {H_*(\{f\ge 0\},\KM)}$). 
We say that a point $x \in J^1(N)$ is {\em above} some subset of the manifold $N$ if its image under the natural projection $J^1(N) \rightarrow N$ belongs to this subset.

\begin{thm} \label{thm4}
For any positive Legendrian manifold  $L \subset (J^1(N), \zeta)$ in general
position with respect to $\Lambda$, 
there exists at least $b(f)$ points of intersection of $L$ with
$\Lambda_+$ lying above $M\setminus \partial M$.

More precisely, for a generic positive Legendrian manifold $L$, there exists at least $b(f)$ 
different {\em positive} numbers $\lambda_1,\dots,\lambda_{b(f)}$ such that
$L$ intersects each manifold $j^1(\lambda_i f)$ above $M\setminus \partial M$.
\end{thm}

\begin{rem}
Theorem \ref{thm4} implies the Morse estimate for the number of
critical points of a Morse function $F$ on $N$. This can be seen as follows.
By adding a sufficiently large constant to $F$, one can assume that $L=j^1(F)$
is a positive Legendrian manifold. 
If $f$ is a constant positive function, then $M=N$, and intersections
of $L=j^1F$ with $\Lambda$ are in one to one correspondence with the
critical points of $F$. Furthermore, $F$ is Morse if and only if $L$
is transversal to $\Lambda$.

\medskip

In fact, one can prove that Theorem \ref{thm4} implies (a weak form of) Arnold's conjecture for
Lagrangian intersection in cotangent bundles, proved by Chekanov \cite{Ch} in its Legendrian version.
This is no accident: our proof rely on the main ingredient of
Chekanov's proof: the technique of generating families (see Theorem \ref{Chgen}).

\medskip

One can also prove that Theorem  \ref{thm4} implies Theorem \ref{thm1}.
Theorem \ref{thm5} below, which in turn implies Theorem \ref{thm2},
is also a direct consequence of Theorem \ref{thm4}.

\end{rem}

\begin{thm}\label{thm5}
Consider a line in $\RM^n$. Denote by $\Lambda$
the union of all the fibers of
$\pi\colon ST^*{\RM}^n \to {\RM}^n $ above this line.
Consider one of these fibers and a positive path starting from this fiber. 
The end of this positive path is a Legendrian sphere. 
This sphere must intersect $\Lambda$ in at least $2$ points.
\end{thm}

\subsection{An application to positive isotopies in homogeneous neighborhoods of a surface in a tight contact 3 manifold}
In Theorem \ref{thm4} and in Theorem \ref{thm5}, we observe the following feature:
The submanifold $\Lambda$ is foliated by Legendrian submanifolds. 
We pick one of them, and we conclude that we cannot disconnect it from $\Lambda$ by a positive contact isotopy.
In dimension 3, our $\Lambda$ is a surface foliated by Legendrian curves (in a non generic way).

\medskip

Recall that generically, a closed oriented surface $S$ contained
in a contact $3$-manifold $(M,\xi )$ is {\it convex}: there exists
a vector field tranversal to $S$ and whose flow preserves $\xi$.
Equivalently, a convex surface admits a {\it homogeneous
neighborhood} $U\simeq S\times \RM$, $S\simeq S\times \{ 0\}$,
where the restriction of $\xi$ is $\RM$-invariant. Given such an
homogeneous neighborhood, we obtain a smooth, canonically
oriented, multicurve $\Gamma_U \subset S$, called the {\it
dividing curve} of $S$, made of the points of $S$ where $\xi$ is
tangent to the $\RM$-direction. It is automatically transversal to
$\xi$. According to Giroux \cite{Gi}, the dividing curve $\Gamma_U$ does not
depend on the choice of $U$ up to an isotopy amongst the
multicurves transversal to $\xi$ in $S$. The {\it characteristic
foliation} $\xi S$ of $S\subset (M,\xi )$ is the integral foliation
of the singular line field $TS\cap \xi$.

\medskip

Let $S$ be a closed oriented surface of genus $g(S)\geq 1$ and
$(U,\xi )$, $U\simeq S\times \RM$, be an homogeneous neighborhood
of $S\simeq S\times \{0\}$. The surface $S$ is $\xi$-convex, and
we denote by $\Gamma_U$ its dividing multicurve. We assume that $\xi$
is tight on $U$, which, after Giroux, is the same than to say that no component
of $\Gamma_U$ is contractible in $S$.

\begin{thm}
\label{thm6}
Assume $L$ is a Legendrian curve in $S$ having minimal geometric
intersection $2k>0$ with $\Gamma_U$. If $(L_s)_{s\in [0,1]}$ is a
positive Legendrian isotopy of $L=L_0$ 
then $\sharp (L_1 \cap S )\geq 2k$.
\end{thm}

\begin{rem}
The positivity assumption is essential: if we push $L$ in the
homogeneous direction, we get an isotopy of Legendrian curves
which becomes instantaneously disjoint from $S$. If $k=0$, this is
a positive isotopy of $L$ that disjoints $L$ from $S$.
\end{rem}

\begin{rem}
For a small
positive isotopy, the result is obvious. Indeed, $L$ is an
integral curve of the characteristic foliation $\xi S$ of $S$,
which contains at least one singularity in each component of
$L \setminus \Gamma_U$. For two consecutive components, the singularities
have opposite signs. Moreover, when one moves $L$ by a small
positive isotopy, the positive singularities are pushed in $S
\times \RM^+$ and the negative ones in $S \times \RM^-$. Between two
singularities of opposite signs, we will get one intersection with
$S$.
\end{rem}

\medskip

The relationship with the preceeding results is given by the following corollary of 
theorem \ref{thm4}, applied with $N=S^1$ and $f(\theta)= \cos (k\theta)$, for some fixed $k \in \NM $.
In this situation, the surface $\Lambda$ of Theorem \ref{thm4} will be called $\Lambda_k$.
It is an infinite cylinder foliated by Legendrian circles. 
Its characteristic foliation $\xi \Lambda_k$ has $2k$ infinite lines of singularities.
The standard contact space $(J^1(S^1), \zeta)$ is itself an homogeneous neighborhood of $\Lambda_k$,
and the corresponding dividing curve consists in $2k$ infinite lines, alternating with the lines of singularities.

\medskip

Let $L_0=j^10 \subset \Lambda_k$. Theorem \ref{thm4} gives:

\begin{cor}
\label{ex:fund}
Let $L_1$ be a generic positive deformation of $L_0$. Then $\sharp\{L_1 \cap \Lambda_k\}\geq 2k$. 
\end{cor}

Indeed, there are $k$ intersections with $\Lambda_{k,+}$, and $k$ other intersections which are obtained in a similar way with the function $-f$.
\qed

\medskip

This corollary will be the building block to prove Theorem \ref{thm6}.


\subsection{Organization of the paper}
This paper is organized as follows.
The proof of Theorem $3$, which in a sense shows that the hypothesis of Theorems \ref{thm1} and \ref{thm2}
are optimal, consists essentially in a collection of explicit constructions.
It is done in the next section (\ref{proof3}) and it might serve as an
introduction to the main notions and objects discussed in this paper.
The rest of the paper is essentially devoted to the proof of Theorems
1, 2, 4, 5 and 6, but contains a few statements which are more general
than the theorems mentioned in this introduction.


\subsection{Acknowledgements}
This work was motivated by a question of Yasha Eliashberg.
This paper is an based on an unpublished preprint of 2006 \cite{CFP}.
Since then, Chernov and Nemirovski proved a statement which generalises our Theorems 1 and 2 \cite{CN1,CN2}, 
and they found new applications of this to causality problems in space-time.  
All this is also related to the work of Bhupal \cite{Bh} and to the
work of Sheila Sandon \cite{Sa}, 
who reproved some results of \cite{EKP} using the generating families
techniques.

\medskip

Vincent Colin is partially supported by the ANR Symplexe, the ANR Floer
power and the Institut universitaire de France.
Petya Pushkar is partially supported by RFBR grant 08-01-00388. 


\section{Proof of theorem \ref{thm3} }\label{proof3}

\subsection{A positive loop}\label{posLoop}
In order to prove statement { \bf i)} of Theorem \ref{thm3},
we begin by the description of a positive loop in the space of Legendrian embeddings in $J^1(S^1)$.
Due to Theorem \ref{thm1}, this cannot happen in the component of the zero section $j^10$.

\medskip

Take $\epsilon >0$ and consider a
Legendrian submanifold $L$ homotopic to $j^10$ and embedded in the half-space $\{p>2\epsilon \}\subset J^1(S^1)$.
Consider the contact flow $\varphi_t\colon (q,p,u) \to (q-t,p,u-t\epsilon), t\in \RM$.
The corresponding contact Hamiltonian
$h(q,p,t)=-\epsilon+p$ is positive near $\varphi_t(L)$, for all $t \in \RM$, 
and hence, $\varphi_t(L)$ is a positive path.

\medskip

On the other hand, one can go from $\varphi_{2\pi}(L)$ back to $L$ just
by increasing the $u$ coordinate, which is also a positive path.
This proves statement {\bf i)} of Theorem \ref{thm3}.

\begin{figure}
\label{fig2}
\begin{center}
\includegraphics{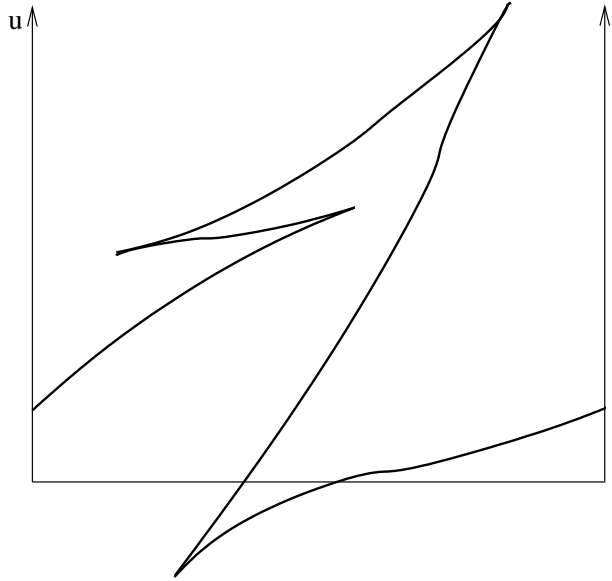}
\caption{The front projection of an $L \subset \{p>2 \epsilon\}$,
  wxhich is homotopic to $j^10$ through Legendrian immersions.}
\end{center}  
\end{figure}

\subsection{} \label{thm3A} We now consider statement {\bf ii)}.
Take $L \subset \{p>2\epsilon \}$ as above,
but assume in addition that $L$ is homotopic to $j^10$ through Legendrian immersions.
Such a $L$ exists (one can show that the $L$ whose front projection is depicted in
fig. \ref{fig2} is such an example), but cannot be Legendrian isotopic to $j^10$,
since, by \cite{Ch}, it would intersect $\{p=0\}$).

\medskip

\paragraph{\bf Step 1.} The homotopy between $j^10$ and $L$ can be transformed into a positive path of Legendrian immersions between $j^10$ and a vertical translate $L'$ of $L$, by combining it with an
 upwards translation with respect to the $u$ coordinate.

\medskip

\paragraph{\bf Step 2.} Then, using the flow $\varphi_t$ (defined in \ref{posLoop}) for $t \in
[0, 2k\pi]$ with $k$ big enough, one can reach reach another translate
$L''$ of $L$, on which the $u$ coordinate can be arbitrarily low.

\medskip

\paragraph{\bf Step 3.} Consider now a path of Legendrian immersions from $L$ to $j^10$. 
It can be modified into a positive path between $L''$ and $j^10$, like in Step~1.

\medskip

This proves statement {\bf ii)} of Theorem \ref{thm3}.

\subsection{} The proof of statement {\bf 0)} uses again the same idea.
Given any compact Legendrian submanifold $L \subset J^1(\RM^n)$, there exists $L'$, which is Legendrian isotopic to $L$ and which is contained into the half-space $p_1>\epsilon>0$, for some system $(p_1,\dots,p_n,q_1,\dots, q_n)$ of canonical coordinates on $T^*(\RM^n)$.
It is possible to find a positive path between $L$ and a sufficiently high vertical translate $L''$ of $L'$.
Because $p_1>\epsilon$, one can now slide down $L''$ by a positive path 
as low as we want with respect to the $u$ coordinate, as above.
Hence we can assume that $L$ is connected by a positive path of embeddings to some $L'''$, 
which is a vertical translate of $L'$, on which the $u$ coordinate is very negative. 
So we can close this path back to $L$ in a positive way. 

\subsection{} We now prove statement {\bf iii)}.
Consider two points $x$ and $y$ on the surface $N$, an embedded path from
$x$ to $y$, and an open neighborhood $U$ of this path, diffeomorphic
to $\RM^2$. Hence it is
enough to consider the particular case $N=\RM^2$. We consider this
case below.

\medskip

\subsubsection{The hodograph transform} \label{hodograph}
We now recall the classical "hodograph" contactomorphism \cite{Acusps} which identifies 
$(ST^*\RM^2, \zeta_1)$ and $(J^1(S^1), \zeta)$, and more generally
$(ST^*\RM^n, \zeta_1)$ and ($J^1(S^{n-1}), \zeta)$. 
The same trick will be used later to prove Theorems 2 and 5 (sections \ref{proof2} and \ref{proof4}).

\medskip

Fix a scalar product $\langle.,.\rangle$ on $\RM^n$ and identify the sphere $S^{n-1}$ with the standard unit sphere in $\RM^n$. 
Identify a covector at a point $q \in S^{n-1}$ with a vector in the hyperplane tangent to the sphere at $q$
(perpendicular to $q$). 
Then to a point $(p,q,u)\in J^1(S^{n-1})=T^*S^{n-1}\times\RM$ we associate the cooriented contact
element at the point $uq+p\in\RM^n$, which is parallel to $T_qS^{n-1}$, and cooriented by $q$.

\medskip

One can check that the fiber of $\pi \colon ST^*\RM^n \rightarrow \RM$ over some point $x \in \RM^n$ is the
image by this contactomorphism of $j^1l_x$, 
where $l_x: S^{n-1} \rightarrow \RM$, $q \mapsto \langle x,q\rangle$.

\subsubsection{End of the proof of Theorem \ref{thm3} {\bf iii)}}
One can assume that $x=0 \subset \RM^2$. 
The case when $x=y$ follows directely from Theorem \ref{thm3} {\bf ii)} via the
contactomorphism described above. The fiber $\pi^{-1}(x)$ corresponds to $j^10$.

\medskip

Suppose now that $x \neq y$. 
We need to find a positive path of Legendrian immersions in
$(J^1(S^1), \zeta)$ between $j^10$ and $j^1l_y$.

\medskip

To achieve this, it is enough to construct a  positive path of
Legendrian immersions between $j^10$ and a translate of $j^10$ that
would be entirely below $j^1l_y$, with respect to the $u$ coordinate.
This can be done as in \ref{thm3A}, just by decreasing even
more the $u$ coordinate like in step 2. 
This finishes the proof of Theorem \ref{thm3}. \qed

\bigskip


\section{Morse theory for generating families quadratic at infinity}

\subsection{Generating families} 
We briefly recall the construction of a generating family for a Legendrian manifold 
(the details can be found in~\cite{AG}).
Let $\rho \colon E \to N$ be a smooth fibration over a smooth manifold $N$, with fiber $W$. 
Let $F\colon E \to \RM$ be a smooth function. 
For a point $q$ in $N$ we consider the set $B_q \subset \rho^{-1}(q)$ 
whose points are the critical points of the restriction of $F$ to the fiber $\rho^{-1}(q)$. 
Denote $B_{F}$ the set $B_F=\bigcup_{q\in N} B_q \subset E$. 
Assume that the rank of the matrix $ (F_{wq},F_{ww})$ ($w$, $q$ are
local coordinates on the fiber and base respectively) 
formed by second derivatives is maximal (that is, equal to the
dimension of~$N$) at each point of $B_{F}$. This condition holds for a
generic~$F$ and does not depend on the choice of the local coordinates
$w,q$. 

\medskip

The set $B_{F}\subset W$ is then a smooth submanifold of the same dimension as~$N$, 
and the restriction of the map $$(q,w)\stackrel {l_{F}}\longermapsto (q,d_{N}(F(q,w)), F(q,w)),$$ 
where $d_N$ denotes the differential along $N$, to $B_{F}$ defines a
Legendrian immersion of $B_{\?F}$ into~$(J^1(M), \zeta)$. 
If this is an embedding (this is generically the case), 
then $F$ is called a {\it generating family} of the Legendrian submanifold $L_F = l_{\?F}(B_{\?F})$.

\medskip

A point $x \in J^1(N)$ is by definition a triple consisting in a point
$q(x)$ in the manifold $N$, a covector $p(x)\in T^*_{q(x)}N$ and a real number $u(x)$. 
A point $x\in L$ will be called a {\it critical point} of the
Legendrian submanifold $L\subset (J^1(N), \zeta)$ if $p(x)=0$. 
The value of the $u$ coordinate at a critical point of a Legendrian
manifold $L$ will be called a {\em critical value} of $L$. The set of
all critical values will be denoted by $Crit(L)$. 

\medskip

Observe that, for a manifold $L=L_F$ given by a generating family $F$ the set 
$Crit(L_F)$ coincides with the set of critical values of the generating family~$F$.
 
\medskip 
 
We call a critical point $x\in L$ {\it nondegenerate} if $L$ intersects the manifold given by the equation
$p=0$ transversally at $x$. 
If an embedded Legendrian submanifold $L_F$ is given by a generating family $F$, then the
non-degenerate critical points of $F$ are in one to one correspondance with the non-degenerate critical points of $L_F$.

\medskip

We describe now the class of generating families we will be working with. 
Pick a closed manifold $E$ which is a fibration over some closed manifold $N$. 
A function $F\colon E\times \RM^K \?\?\? \to \RM$ is called {\it
  $E$-quadratic at infinity\/} if it is a sum of a non-degenerate
quadratic form $Q$ on $\RM^K$ and a function on $E\times\RM^K$ with
bounded differential (i.e. the norm of the differential is uniformly bounded
for some Riemannian metric which is a product of a Riemannian metric
on $E$ and the Euclidean metric on $\RM^K$). 
This definition does not depend on the choice of the metrics.
If a function which is $E$-quadratic at infinity is a generating
family (with respect to the fibration $E\times \RM^K \to N$), then we
call it a {\em generating family $E$-quadratic at infinity}.

\subsection{Morse theory for generating families $E$-quadratic at infinity} 
We gather here some results from Morse theory which will be needed later. 
Let $E\to N$ be a fibration, $E$ is a closed manifold. Consider a function $F$, $E$-quadratic at infinity.
Denote by $F^a$ the set $\{F\le a\}$. 
For sufficiently big positive numbers $C_1<C_2$, the set $F^{-C_2}$ is a 
deformation retract of  $F^{-C_1}$. 
Hence the homology groups $H_*(F^a, F^{-C}, \KM)$ depend only on $a$.
We will denote them by $H_*(F,a)$. 
It is known (see \cite{CZ}) that for sufficiently big $a$, $H_*(F,a)$ is isomorphic to $H_*(E, \KM)$.

\medskip

For any function $F$ which is $E$-quadratic at infinity, and any integer\\ $k\in \{1,\dots,\dim H_*(E,\KM)\}$, 
we define a {\em Viterbo number}
$c_k(F)$ by
$$
c_k(F)=\inf\{c|\dim i_*(H_*(F,c))\ge k\},
$$
where $i_*$ is the map induced by the natural inclusion $F^c\to F^a$,
when $a$ is a sufficiently big number. Our definition is similar to
Viterbo's construction \cite{Vi} in the symplectic setting. The
following proposition is an adaptation of \cite{Vi}:

\begin{prop}
\label{Morseprop}
{\bf i.}~Each number $c_k(F)$, $k\in\{1,\dots,\dim H_*(E, \KM)\}$ is a critical value of $F$,
and if $F$ is an excellent Morse function (i.e all its critical points are non-degenerate and all critical values are
different) then the numbers $c_k(F)$ are different.

\medskip

{\bf ii.}~Consider a family $F_{t, t\in [a,b]}$ of functions which are all $E$-quadratic at infinity. 
For any $k\in\{1,\dots,\dim H_*(E, \KM)\}$ the number $c_k(F_t)$ depends on $t$ continuously. 
If the family $F_{t, t\in [a,b]}$ is generic (i.e. intersects the discriminant
consisting of non excellent Morse functions transversally at its smooth points) 
then $c_k(F_t)$ is a continuous piecewise smooth function with a finite number of singular points. \qed

\end{prop}

\begin{rem}
At this moment, it is unknown wether $c_i(L_F)$ depends on $F$ for a given Legendrian manifold $L=L_F$.
Conjecturally there should be a definition of some analogue of $c_i$ in terms of augmentations on relative contact homology.
\end{rem}


\section{Proof of Theorem \ref{thm1}}

We will in fact prove Theorem \ref{thm7} below, which is more general than Theorem \ref{thm1}. 
Fix a closed (compact, without boundary) manifold $N$ and a smooth fibration $E\to N$ such that $E$ is compact. 
A Legendrian manifold $L\subset (J^1(N), \zeta)$ will be called a {\em $E$-quasifunction} if it is Legendrian isotopic to a
manifold given by some generating family $E$-quadratic at infinity. 
We say that a connected component $\mathcal{L}$ of the
space of Legendrian submanifolds in $(J^1(N), \zeta)$ is {\em $E$-quasifunctional} if $\mathcal{L}$ contains an
$E$-quasifunction. 
For example, the component $\mathcal{L}$ containing the one jets
extensions of the smooth functions on $M$ is $E$-quasifunctional, with
$E$ coinciding with $N$ (the fiber is just a point). 

\begin{thm} \label{thm7}
\label{noqloops} 
An $E$-qua\-si\-func\-ti\-onal component contains no closed positive path.
\end{thm}

\subsection{} 
The proof of Theorem 7 will be given in \ref{pfthm7}. 
It will use the following generalization of Chekanov's theorem (see \cite{P1}), and proposition \ref{prop} below.

\begin{thm}
\label{Chgen} Consider a Legendrian isotopy $L_{t, t\in [0,1]}$ such
that $L_0$ is an $E$-quasifunction.
Then there exist a number $K$ and a smooth family of
functions $E$-quadratic at infinity $F_t\colon E\times\RM^K \to \RM$, such that for any $t\in[0,1]$, $F_t$ is a generating family of $L_t$. \qed
\end{thm}

\medskip

Note that it follows from Theorem~\ref{Chgen} that {\em any} Legendrian manifold in some $E$-quasifunctional component is in fact an $E$-quasifunction.

\medskip

Consider a positive path $L_{t,t\in [0,1]}$ given by a family $F_{t,t\in[0,1]}$ of E-quadratic at infinity generating families. 
We are going to prove the following inequality:

\medskip

\begin{prop}\label{prop} The Viterbo numbers of the family $F_t$ are monotone increasing functions with respect to $t$:
$c_{i}(F_0)<c_{i}(F_1)$ for any $i\in\{1,...,\dim H_*(E)\}$.
\end{prop}

\subsection{Proof of Proposition \ref{prop}}
Assume that the inequality is proved for a generic family. 
This, together with continuity of Viterbo numbers, gives us a weak inequality $c_{i,M}(F_0)\le c_{i,M}(F_1)$, 
for any family. 
But positivity is a $C^{\infty}-$open condition, so we can perturb the initial 
family $F_t$ into some family $\widetilde{F}_t$ coinciding with $F_t$ when $t$ is sufficiently close to $0,1$, such that $\widetilde{F}_t$ still
generates a positive path of legendrian manifolds and such that the family $\widetilde{F}_{t, t \in [1/3,2/3]}$ is
generic. We have
$$
c_{i}(F_0)=c_{i,M}(\widetilde{F}_0)\le c_{i}(\widetilde{F}_{1/3}) < c_{i}(\widetilde{F}_{2/3})\le
c_{i}(\widetilde{F}_1)=c_{i}(F_1),
$$
and hence inequality is strong for all families.

\medskip

We now prove the inequality for generic families.
Excellent Morse functions form an open dense set in the space of all E-quadratic at infinity functions on $N\times \RM^K$. 
The complement of the set of excellent Morse functions forms a discriminant, which is a singular hypersurface. 
A generic one-parameter family of E-quadratic at infinity functions $F_t$ on $E \times \RM^K$ 
has only a finite number of transverse intersections with the discriminant in its smooth points,
and for every $t$ except possibly finitely many, the Hessian $d_{ww}F_{t}$ is non-degenerate at every critical point of the function $F_t$.

\medskip

We will use the notion of {\em Cerf diagram} of a family of functions $g_{t,t\in[a,b]}$ on a smooth manifold. 
The Cerf diagram is a subset in  $[a,b]\times \RM$ consisting of all the pairs of type $(t,z)$, where $z$ is a critical value of $g_t$. 
In the case of a generic family of functions on a closed manifold, the
Cerf diagram is a curve with non-vertical tangents 
everywhere, with a finite number of transversal self-intersections and cuspidal points as singularities. 

\medskip

The graph of the Viterbo number $c_{i}(F_t)$ is a subset of the Cerf diagram of the family $F_t$. 
To prove the monotonicity of the Viterbo numbers, it is sufficient to
show that the Cerf diagram of $F_t$ has a positive slope at every
point except finite set.
The rest of the proof of Proposition \ref{prop} is devoted to that.

\medskip

We say that a point $x$ on a Legendrian manifold $L\subset J^1(N)$ is non-vertical if the differential of the
natural projection $L \to N$ is non-degenerate at $x$. 
Let $L_t$ be a smooth family of Legendrian manifolds and $x(t_0)=(p(t_0), q(t_0),u(t_0))$ a non-vertical point. 
By the implicit function theorem, there exists a unique family
$x(t)=(p(t),q(t),u(t))$, defined for $t$ sufficiently 
close to $t_0$, such that $x(t)\in L_t$ and $q(t)=q(t_0)$. 
We call the number $\frac {d}{dt}\big|_{t=t_0}u(t)$ {\em vertical speed} of the point $x(t_0)$.

\begin{lem} \label{l1}
For a positive path of Legendrian manifolds, the vertical speed of every non-vertical point is positive. \qed
\end{lem}

Consider a path $L_t$ in the space of legendrian manifolds given by a ge\-ne\-rating family $F_t$. 
Consider the point $x(t_0)\in L_{t_0}$ and the point $(q,w)\in N\times \RM^{K}$ such that 

$$d_wF_{t_0}(q,w)=0, x(t_0)=(p,q,u), p=d_qF_{t_0}(q,w), u=F_{t_0}(q,w).$$

Then $x$ is non-vertical if and only if the hessian
$d_{ww}F_{t_0}(q,w)$ is non-degenerate. For such a point $x$, the following lemma holds:

\begin{lem} \label{l2}
The vertical speed at $x$ is equal to $\frac
{d}{dt}\big|_{t=t_0}F_t(q,w)$ \qed.
\end{lem}

Let $G_t$ be a family of smooth functions and assume that the point $z(t_0)$ is a Morse critical point for $G_{t_0}$. 
By the implicit function theorem, for each $t$ sufficiently close to $t_0$, 
the function $G_t$ has a unique critical point $z(t)$ close to $z_0$, and $z(t)$ is a smooth path. 

\begin{lem} \label{l3}
The speed of the critical value $\frac {d}{dt}\big|_{t=0}G_t(z(t))$ is equal to $\frac
{d}{dt}\big|_{t=0}G_t(z(t_0))$.
\end{lem}

Indeed,  $\frac {d}{dt}\big|_{t=t_0}G_t(z(t)) = \frac{d}{dt}\big|_{t=t_0}G_t(z(t_0))+\frac
{\partial G_t}{\partial z}(z(t_0)) \cdot \frac {dz}{dt}(t_0)$. \qed

\medskip

At almost every point on the Cerf diagram, the slope of the Cerf diagram at this point is the speed of
a critical value of the function $F_t$. By lemma \ref{l3} and lemma \ref{l2}, it is the
vertical speed at some non-vertical point. By lemma \ref{l1} it is
positive. This finishes the proof of Proposition \ref{prop}. \qed

\medskip

\subsection{Proof of Theorem \ref{thm7}} \label{pfthm7}
Suppose now that there is a closed positive loop $L_{t, t\in [0,1]}$ 
in some $E$-qua\-si\-func\-ti\-onal component $\mathcal{L}$.
The condition of positivity is open. 
We slightly perturb the loop $L_{t, t\in [0,1]} $ 
such that $Crit(L_0)$ is a finite set of cardinality $A$.
Note that $A>0$, since $L_0$ is a $E$-quasifunction. 
Consider the $A$-th multiple of the loop $L_{t, t\in [0,1]}$. 
By Theorem~\ref{Chgen}, $L_t$ has a generating family $F_t$,
for all $t\in [0,A]$. 
By Proposition~\ref{prop}, we have that
$$c_1(\widetilde{F_0})<c_1(\widetilde{F_1})<\dots<c_1(\widetilde{F_A}).$$
All these $A+1$ numbers belong to the set $Crit(L_0)$.
This is impossible due to the cardinality of this set. 
This finishes the proof of Theorem~\ref{thm7} and hence of Theorem
$\ref{thm1}$. \qed

\medskip



\subsection{Proof of Theorem \ref{thm2}} \label{proof2}
Theorem \ref{thm2} is a corollary of Theorem \ref{thm1}, 
via the the contactomorphism between $(ST^*(\RM^n), \zeta_1)$ and $(J^1(S^{n-1}), \zeta)$ we have seen in \ref{hodograph}.

\medskip

Consider the fiber $\pi^{-1}(x)$ of the fibration $\pi\colon ST^*\RM^n \to \RM^n$. 
It corresponds to a Legendrian manifold $j^1l_x\subset J^1(S^{n-1})$,
where $l_x$ is the function $l_x=\langle q,x\rangle$. 
It is a Morse function for $x\ne 0$, and has only two critical points and two critical values $\pm||x||$. 
The critical points of $l_x$ are non-degenerate if $x\ne 0$. 
It follows from Proposition~\ref{Morseprop} that $c_1(F)=-||x||$, $c_2(F)=||x||$. 
Indeed, any small generic Morse perturbation of $F$ has two critical
points with critical values close to $\pm||x||$. Viterbo numbers for
this perturbation should be different. Hence by continuity
$c_1(F)=-||x||$, $c_2(F)=||x||$. 

Viterbo numbers for $j^1l_0$ are equal to zero, because $Crit(S(0))=\{0\}$. 
The existence of a positive path would contradict the monotonicity
(Proposition \ref{prop}) of Viterbo numbers. \qed


\section{Morse theory for positive Legendrian submanifolds}
\label{proof4}
In this section we prove Theorem \ref{thm4} and deduce Theorem \ref{thm5} from it. 
We need first to generalize some of the previous constructions and results to the case of manifolds with boundary.

\medskip

Let $N$ be a compact closed manifold. 
Fix a function $f\colon W\to\RM$ such that $0$ is a regular value of $f$. 
Denote by $M$ the set $f^{-1}([0,+\infty[)=\{f\ge0\}$. 
Denote by $b(f)=\dim_{\KM} {H_*(M)}$ the dimension $H_*(M)$ (all the homologies here and below are
counted with coefficients in a fixed field $\KM$).

\subsection{Viterbo numbers for manifolds with boundary}
The definition of the Viterbo numbers for a function quadratic at infinity on a manifold 
with boundary is the same as in the case of closed manifold. 
We repeat it briefly. Given a function $F$ which is quadratic at infinity, we define the Viterbo
numbers $c_{1,M}(F),...,c_{b(f),M}(F)$ as follows.

\medskip

A {\em generalized critical value} of $F$ is a real number which is a critical for $F$ or for
the restriction $F|_{\partial M\times \RM^N}$. Denote by $F^a$ the set
$\{(q,w)|F(q,w) \leq a \}$. 
The homotopy type of the set $F^a$ is changed only if $a$ passes through a generalized critical value. 
One can show that, for sufficiently big $K_1,K_2>0$, the homology
of the pair $(F^{K_1},F^{-K_2})$ is independent of $K_1,K_2$, and naturally isomorphic, by the Thom isomorphism,
to $H_{*-\ind Q}(M)$.
So, for any $a \in \RM$ and sufficiently large $K_2$ the projection
$H_*(F^a,F^{-K_2}) \to H_{*-\ind Q}(M)$ 
is well defined and independent of $K_2$. 
Denote the image of this projection by~$I(a)$.

\begin{defi} The Viterbo numbers are
$$c_{k,M} (F) = \inf\{c| \dim I(c) \ge k \}, k \in \{1, \dots, b(f))\}.$$
\end{defi}

Any Viterbo number $c_{k,M}(F)$ is a generalized critical value of the function $F$.
Obviously, $c_{1,M}(F)\le ...\le c_{b(f),M}(F)$. 
For any continuous family $F_t$ of quadratic at infinity
functions, $c_{i,M}(F_t)$ depends continuously on $t$. 

\subsection{Proof of Theorem \ref{thm4}}
Consider a $1$-parameter family of quadratic at infinity functions $F_{t, t\in[a,b]}\colon N\times\RM^{K}\to
\RM$, such that $F_t$ is a generating family 
for the Legendrian manifold $L_t$ and such that the path $L_{t,t\in [a,b]}$ is positive. 
We will consider the restriction of the function $F_t$ to $M\times\RM^{K}$ and denote it by $F_t$ also.
The following proposition generalizes Proposition~\ref{prop}.

\begin{prop}
\label{prop1} 
The Viterbo numbers of the family $F_t$ are monotone increasing:
$c_{i,M}(F_a)<c_{i,M}(F_b)$ for any $i\in\{1,...,b(f)\}$.
\end{prop}

The difference with Proposition \ref{prop} is that the Cerf
diagram of a generic family has one more possible singularity. 
This singularity corresponds to the case when a Morse critical point
meets the boundary of the manifold. 
In this case, the Cerf diagram is locally diffeomorphic to a
parabola with a tangent half-line. \qed

\medskip

We now prove theorem \ref{thm4}. 
Consider the 1-parameter family of functions $H_{\lambda\ge 0}$, $$H_\lambda(q,w)=F_1(q,w)-\lambda f(q)$$ 
on $M\times \RM^K$. The manifold $L$ intersects $j^1\lambda_0f$ at some point above $M$ if and only if the
function $H_{\lambda_0}$ has $0$ as an ordinary critical value (not a critical value of
the restriction to the boundary).

\medskip

Consider the numbers $c_{k,M}(H_\lambda)$. 
By Proposition \ref{prop1}, $c_{k,M}(H_0)>0$. 
For a sufficiently big value of $\lambda$, each of them is negative.
To show that, consider a sufficiently small $\varepsilon >0$ belonging to the component of the regular values of $f$ which contains $0$. 
Denote by $M_1\subset M$ the set $\{f\ge\varepsilon\}$. 
The manifold $M_1$ is diffeomorphic to the manifold $M$, and the inclusion map is an homotopy equivalence. 
Denote by $G_\lambda$ the restriction of $H_\lambda$ to the $M_1\times\RM^K$.
Consider the following commutative diagram:
$$
\begin{CD}
H_*(G_\lambda^a,G_\lambda^{-K_2})@>i_1>>H_*(G_\lambda^{K_1},G_\lambda^{-K_2})@>Th_1>>H_{*-\ind(Q)}(M_1)\\
@Vj_1VV @. @Vj_2VV\\
H_*(H_\lambda^a,H_\lambda^{-K_2})@>i_2>>H_*(H_\lambda^{K_1},H_\lambda^{-K_2})@>Th_2>>H_{*-\ind(Q)}(M)
\end{CD}
$$
where $K_1,K_2$ are sufficiently big numbers, $Th_1,Th_2$ denote Thom isomorphisms and $i_1,i_2,j_1,j_2$ are the maps induced by the natural inclusions. 
It follows from the commutativity of the diagram and from the fact that $j_2$ is an isomorphism that
$c_{k,M_1}(G_\lambda)\ge c_{k,M}(H_\lambda)$ for every~$k$.

\medskip

For sufficiently big $\lambda$ and for every $q\in M_1$, 
the critical values of the function $G_\lambda$ restricted to $q\times \RM^K$ are negative.
Hence all generalized critical values of $G_\lambda$ are negative.
It follows that all the numbers $c_{k,M_1}(G_\lambda)$ are negative, and the same holds for $c_{k,M}(H_\lambda)$. 
We fix $\lambda_0$ such that $c_{k,M}(H_{\lambda_0})<0$ for every $k
\in \{1,\dots, b(f) \}$.

\medskip

Consider now  $c_{k,M}(H_\lambda)$ as a function of $\lambda\in[0,\lambda_0]$.
We are going to show that its zeroes correspond to the intersections above $M\setminus \partial M$.
For a manifold $L_1$ in general position,  all the generalized critical values of $F_1$ are non-zero. 
In particular all the critical values of the function the $F_1|_{\partial M\times \RM^K}$ are non-zero. 
The function $F_1|_{\partial M\times \RM^K}$ coincides with ${H_\lambda |} _{{\partial M\times \RM^K}}$ since $f=0$ on $\partial M$. 
Hence, if zero is a critical value for $H_\lambda$, then it is an ordinary critical value at some inner point. 
This finishes the proof of Theorem \ref{thm4}. \qed

\medskip

\begin{rem}The function $c_{i,M}(H_{\lambda})$ can be constant
  on some sub-intervals in $]0,\lambda_0[$, even for a generic function
  $F_1$. Indeed, the critical values of the restriction of $H_\lambda$
  to $\partial M\times\RM^K$ do not depend on $\lambda$. It is
  possible that $c_{i,M}(H_\lambda)$ is equal to such a critical value
  for some $\lambda$'s.
\end{rem}

The following proposition concerns the case of a general (non
necessarily generic) positive Legendrian manifold. 
We suppose again that $f$ is a function having $0$ as regular value and that $L$ is a positive manifold.

\begin{prop} \label{common}
For any connected component of the set $M=\{f\ge 0\}$ there exists a 
positive $\lambda$ such that $L$ intersects with $j^1\lambda f$ above this component.
\end{prop}

Consider a connected component $M_0$ of the
manifold $M$. It is possible to replace $f$ by some function
$\tilde{f}$ such that $0$ is a regular value for $\tilde{f}$,
$\tilde{f}$ coincides with $f$ on $M_0$ and $\tilde{f}$ is negative on
$N\setminus M_0$. 
We consider $c_{1,M_0}(F_1-\lambda \tilde{f})$ as a function of
$\lambda$. 
It is a continuous function, positive in some neighborhood of zero,
and negative for the big values of $\lambda$. 

\medskip

Fix some $\alpha$ and $\beta$ such that  $c_{1,M_0}(F_1-\alpha
\tilde{f})>0$ and  $c_{1,M_0}(F_1-\beta \tilde{f})<0$.
Assume that for any $\lambda \in [\alpha, \beta]$, $L$ does not
intersect $j^1\lambda\tilde{f}$ above $M_0$. 
Then this is also true for any small enough generic perturbation
$L'$of $L$. Denote by $F'$ a generating family for $L'$. 
Each zero $\lambda_0$ of $c_{1,M_0}(F'-\lambda_0 \tilde{f})$ corresponds to an
intersection of $L'$ with $j^1\lambda_0\tilde{f}$ above $M_0$.
Such a $\lambda_0$ exists by Theorem \ref{thm4}. This is a
contradiction. \qed


\subsection{Proof of Theorem \ref{thm5}} We can suppose that the
origin of $\RM^n$ belongs to the line considered in the statement of
Theorem~\ref{thm5}. Consider now again the contactomorphism of \ref{hodograph}
$(J^1(S^{n-1}, \zeta))=(ST^*(\RM^n), \zeta_1)$. 

\medskip

For such a choice of the origin, the union of all the fibers above
the points on the line forms a manifold of type $\Lambda(f)$, where
$f$ is the restriction of linear function to the sphere $S^{n-1}$. 

\medskip

The manifold $M=\{f\ge 0\}$ has one connected component (it is an hemisphere). 
By Proposition \ref{common} there is at least one intersection of the considered positive Legendrian sphere with $\Lambda_+(f)$. 
Another point of intersection comes from $\Lambda_+(-f)$. 
These two points are different because $\Lambda_+(-f)$ does not intersect with $\Lambda_+(f)$ 


\section{Positive isotopies in homogeneous neighborhoods}


The strategy for proving Theorem~\ref{thm6} is to link
the general case to the case of $\Lambda_k \subset (J^1(S^1) ,\zeta)$. 

\medskip


Let $d=\sharp (L_1 \cap S)$. We first consider the infinite cyclic cover
$\overline{S}$ of $S$ associated with $[L]\in \pi_1 (S)$. The surface
$\overline{S}$ is an infinite cylinder. We call $\overline{U}$ the
corresponding cover of $U$ endowed  with the pullback $\overline{\xi}$
of $\xi$. By construction, $\overline{U}$ is $\overline{\xi}$-homogeneous.
We also call $\overline{L}_s$ a continuous compact lift of $L_s$ in
$\overline{U}$.

\medskip

By compacity of the family $(\overline{L}_s )_{s\in [0,1]}$, we can
find a large compact cylinder $\overline{C} \subset \overline{S}$ such that
for all $s\in [0,1]$, $\overline{L}_s \subset int (\overline{C} \times \RM )$. We
also assume that $\partial \overline{C} \pitchfork \Gamma_{\overline{U}}$.

\medskip

The following lemma shows that in adition we can assume that the boundary of $\overline{C}$ is Legendrian.

\begin{lem} \label{lemmaA} If we denote by $\pi :\overline{S} \times \RM \rightarrow
\overline{S}$ the projection forgetting the $\RM$-factor, we can find a
lift $\overline{C}_0$ of a $C^0$-small deformation of $\overline{C}$ in
$\overline{S}$ which contains $\overline{L}_0$, whose geometric intersection
with $\overline{L}_1$ is $d$ and whose boundary is Legendrian.
\end{lem}


To prove this, we only have to find a Legendrian lift $\gamma$ of a small deformation
of $\partial \overline{C}$, and to make a suitable slide  of $\overline{C}$
near its boundary along the $\RM$-factor to connect $\gamma$ to a
small retraction of $\overline{C} \times \{ 0\}$. The plane field
$\overline{\xi}$ defines a connection for the fibration $\pi :\overline{S}
\times \RM \rightarrow \overline{S}$ outside any small neigbourhood
$N(\Gamma_{\overline{S}} )$ of $\Gamma_{\overline{S}}$. We thus can pick any
$\overline{\xi }$-horizontal lift of $\partial \overline{C}
-N(\Gamma_{\overline{S}} )$. 

\medskip

We still have to connect the endpoints of
these Legendrian arcs in $N(\Gamma_{\overline{S}} )\times \RM$. These
endpoints lie at different $\RM$-coordinates, however this is
possible to adjust since $\overline{\xi}$ is almost vertical in
$N(\Gamma_{\overline{S}} )\times \RM$ (and vertical along
$\Gamma_{\overline{S}} \times \RM$). To make it more precise, we first
slightly modify $\overline{C}$ so that $\partial \overline{C}$ is tangent to
$\overline{\xi} \overline{S}$ near $\Gamma_{\overline{S}}$. Let $\delta$ be  the
metric closure of a component of $\partial \overline{C} \setminus
\Gamma_{\overline{S}}$ contained in the metric closure $R$ of a
component of $\overline{S} \setminus \Gamma_{\overline{S}}$. On $int(R)
\times \RM$, the contact structure $\overline{\xi}$ is given by an
equation of the form $dz +\beta$ where $z$ denotes the
$\RM$-coordinate and $\beta$ is a $1$-form on $int(R)$, such that
$d\beta$ is an area form that goes to $+ \infty$ as we approach
$\partial R$. Now, let $\delta'$ be another arc properly embedded
in $R$ and  which coincides with $\delta$ near its endpoints. If
we take two lifts of $\delta$ and $\delta'$ by $\pi$ starting at the
same point (these two lifts are compact curves, since they
coincide with the characteristic foliation near their endpoints,
and thus lift to horizontal curves near $\Gamma_{\overline{S}}$ where
$\beta$ goes to infinity), the difference of altitude between the
lifts of the two terminal points is given by the area enclosed
between $\delta$ and $\delta'$, measured with $d\beta$. As
$d\beta$ is going to infinity near $\partial \delta =\partial
\delta'$, taking $\delta'$ to be a small deformation of $\delta$
sufficiently close to $\partial \delta$, we can give this
difference any value we want. This proves Lemma \ref{lemmaA} \qed

\medskip

Let $\overline{U}_0 =\overline{C}_0 \times \RM$.

\begin{lem}\label{lemma: embedd} There exists a embedding of $(\overline{U}_0 ,\overline{\xi}
,\overline{L}_0 )$ in $(J^1(S^1) ,\zeta ,\Lambda_k )$ such that the image of
$\overline{L}_1$ intersects $p$ times $A$.
\end{lem}

The surface $\overline{C}_0$ is $\overline{\xi}$-convex and
its dividing set has exactly $2k$ components going from one
boundary curve to the other. All the other components of
$\Gamma_{\overline{C}_0}$ are boundary parallel. Moreover, the curve
$\overline{L}_0$ intersects by assumption exactly once every non
boundary parallel component and avoids the others. Then one can
easily embed $\overline{C}_0$ in a larger annulus $\overline{C}_1$ and
extend the system of arcs $\Gamma_{\overline{C}_0} (\overline{\xi} )$ outside of $\overline{C}_0$ 
by gluing small arcs, in order to obtain a system $\Gamma$ of $2k$ non boundary parallel arcs on $\overline{C}_1$
Simultaneously, we extend the contact structure $\overline{\xi}$ from $\overline{U}_0$,
considered as an homogeneous neighborhood of $\overline{C}_0$, to a
neigbourhood $\overline{U}_1 \simeq \overline{C}_1 \times \RM$ of $\overline{C}_1$.
To achieve this one only has to extend the characteristic foliation, 
in a way compatible with $\Gamma$, and such that the boundary of $\overline{C}_1$ is also Legendrian.
Note that the $\RM$-factor is not changed above $\overline{C}_0$.

\medskip

To summarize, 
$\overline{U}_1$ is an homogeneous neigborhood of $\overline{C}_1$ for the
extension $\overline{\xi}_1$, and $\overline{C}_1$ has Legendrian boundary
with dividing curve $\Gamma_{\overline{C}_1} (\overline{\xi}_1 )=\Gamma$. By
genericity, we can assume that the characteristic foliation of
$\overline{C}_1$ is Morse-Smale. Then, using Giroux's realization
lemma \cite{Gi}, one can perform a $C^0$-small modification of $\overline{C}_1$
relative to $\overline{L}_0 \cup \partial \overline{C}_1$, leading to a
surface $\overline{C}_2$, through annuli transversal to the
$\RM$-direction, and whose support is contained in an arbitrary
small neighborhood of saddle separatrices of
$\overline{\xi}_1 \overline{C}_1$, so that the characteristic foliation of
$\overline{C}_2$ for $\overline{\xi}_1$ is conjugated to $\zeta \Lambda_k$. 
If this support is small enough and if we are in the generic case (which
can always been achieved) where $\overline{L_1}$ doesn't meet the
separatrices of singularities of $\overline{\xi}_1 \overline{C}_1$, 
we get that $\sharp (\overline{L}_1 \cap \overline{C}_2)=\sharp (\overline{L}_1 \cap \overline{C}_1 )=d$. 
As we are dealing with
homogeneous neighborhoods, we see that $(\overline{U}_1 ,\overline{\xi}_1 ,
\overline{L}_0 )$ is conjugated with $(J^1(S^1) ,\zeta ,\Lambda_k )$.
This proves Lamme \ref{lemma: embedd}. \qed

\medskip

The combination of Lemma~\ref{lemma: embedd} and
corollary~\ref{ex:fund} ends the proof of
Theorem~\ref{thm6} by showing that $d \geq 2k$. \qed

\medskip

When $S$ is a sphere the conclusion of
theorem~\ref{thm6} also holds since we are in the
situation where $k=0$. However in this case, we have a more
precise disjunction result.

\begin{thm} Let $(U,\xi )$ be a $\xi$-homogeneous
neighborhood of a sphere $S$. If $\xi$ is tight (i.e. $\Gamma_U$
is connected), then any legendrian curve $L\subset S$ can be made
disjoint from $S$ by a positive isotopy.
\end{thm}


Consider $\RM^3$ with coordinates $(x,y,z)$ endowed
with the contact structure $\zeta =\ker (dz+xdy)$. The radial
vector field $$R=2z\frac{\partial}{\partial z}
+x\frac{\partial}{\partial x}+y\frac{\partial}{\partial y}$$ is
contact. Due to Giroux's realization lemma, the germ of $\xi$ near
$S$ is isomorphic to the germ given by $\zeta$ near a sphere $S_0$
transversal to $R$. Let $L_0$ be the image of $L$ in $S_0$ by this
map. By genericity, we can assume that $L_0$ avoids the vertical
axis $\{ x=0,z=0\}$. Now, if we push $L_0$ enough by the flow of
$\frac{\partial}{\partial z}$, we have a positive isotopy of $L_0$
whose endpoint $L_1$ avoids $S_0$. This isotopy takes place in a
$\zeta$-homogeneous collar containing $S_0$ and obtained by
flowing back and forth $S_0$ by the flow of $R$. This collar
embeds in $U$ by an embedding sending $S_0$ to $S$ and the
$R$-direction to the $\RM$-direction. \qed



\bigskip

\bigskip

\bigskip

{\small

Vincent Colin,
Universit\'e de Nantes, Laboratoire de math\'ematiques Jean Leray, UMR 6629 du CNRS.
email: Vincent.Colin@univ-nantes.fr

\medskip

Emmanuel Ferrand, Universit\'e Pierre et Marie Curie, Institut Math\'ematique de Jussieu, UMR 7586 du CNRS.
email: emmanuel.ferrand@upmc.fr

\medskip

Petya Pushkar, D\'epartement de Math\'ematiques,
Universit\'e Libre de Bruxelles. email: ppushkar@ulb.ac.be

}

\end{document}